\documentstyle[12pt,epsfig]{article}

\newtheorem{thm}{Theorem}
\newtheorem{lem}{Lemma}
\newcommand{\proof}{\vspace*{0.1cm}{\bf Proof.}\hspace{0.2cm}}
\newcommand{\qed}{\hfill Q.E.D.\\ }

\begin{document}
\begin{sloppypar}

\title{
 A Note on the Polytope of Bipartite TSP}
\author{Gergely Kov\'acs$^{a}$,
        Zsolt Tuza$^{b}$
    \thanks{This research was
    supported in part by the Hungarian State and the European Union
    under the grant TAMOP-4.2.2.A-11/1/KONV-2012-0072 },
        B\'ela Vizv\'ari$^{c}$\thanks{Corresponding author, Dept. of IE, 
Eastern Mediterranean University, Famagusta, Turkey, Mersin-10}, Hajieh K. Jabbari$^{c}$}

\maketitle

\parindent 0mm

\small{
$^{a}$ Edutus College, Tatab\'anya, Hungary, kovacs.gergely@edutus.hu\\
$^{b}$ Alfr\'ed R\'enyi Institute of Mathematics,
     Hungarian Academy of Sciences, Budapest, Hungary; and
        Department of Computer Science and Systems
   Technology, University of Pannonia, Veszpr\'em, Hungary, tuza@dcs.uni-pannon.hu\\
$^{c}$ Department of Industrial Engineering, Eastern Mediterranean University, Turkey\\ bela.vizvari@emu.edu.tr, and hajieh.jabbari@cc.emu.edu.tr}

 \vspace{5mm}

\centerline{\bf Abstract}

\parskip 3mm
\parindent 10mm

The main result of this paper is that the polytope of the
bipartite TSP is significantly different from that of the general TSP.
Comb inequalities are known as facet defining ones in the general case.
In the bipartite case, however, many of them are satisfied whenever all degree and
subtour elimination constraints are satisfied, {\em i.e.}\ these comb inequalities
are not facet defining. The inequalities in question belong to the cases where
vertices of one of the two classes occur in less than the half of the intersections
of the teeth and the hand. Such side conditions are necessary, as
 simple example shows that the comb inequality can be violated
when each class has vertices in more than the half of the intersections.

\parskip 2mm


\section{Introduction}
The Traveling Salesman Problem (TSP) has plenty of applications in
operations research. It also has a nice structure and is one of the
most investigated NP-complete problems of combinatorial optimization;
see {\em e.g.}\ the monographs \cite{chvatal}, \cite{gutin}.

Especially a lot of efforts have been made to explore the structure
of the facets of its polytope. The most important facet defining
inequalities are the degree constraints and the subtour elimination
inequalities. The famous Dantzig-Fulkerson-Johnson (DFJ) model
\cite{dfj}  contains only these constraints, and only the characteristic
vectors of the complete
tours satisfy the constraints among the integer
vectors. However, it is well-known that the polytope of the DFJ
model is strictly larger than the TSP polytope, {\em i.e.}\ the convex
hull of the characteristic vectors of the tours.
The next important class of facet defining inequalities is the
set of the
so-called comb inequalities. This type of constraints has several
generalizations including clique-tree, domino, star, path, ladder
and bipartite inequalities. Unfortunately,
a complete description of the TSP polytope is not known \cite{n02}.
Moreover, all of these results concern complete graphs.
However, if the underlying graph is not $K_n$ then the TSP
polytope can be different.

Facet defining inequalities have great algorithmic importance. Branch
and Bound, Branch and Cut, and Cut and Branch methods are used as
general approach to solve TSP. All of them are based on the linear
programming relaxation of TSP. If an optimal solution of the
relaxation is not the characteristic vector of a tour then the bound
can be improved by introducing a violated facet defining inequality.
There are effective methods to find violated subtour elimination
constraints depending on the size of the problem \cite{jv},
\cite{karger},
and \cite{nt02}. There are also some methods to find violated comb
inequalities.

The underlying graph structure affects the difficulty of the particular
TSP. An interesting problem class is bipartite TSP, which
 has application in controlling assembly robots. The
robot moves from cells where parts are stored to an assembly point to
fix the part and goes to another cell after fixing. Thus, its moves have
a bipartite nature. The productivity depends on the length of the route
of the robot arm, and therefore it has to be minimized.

The first paper on
bipartite TSP is \cite{frank}, followed only by rare further publications
\cite{baltz1, baltz2, baltz3}.
%
%
%
%
These papers concentrate on heuristic methods and
 do not discuss the issue of the TSP polytope.
To the best of our knowledge,
 the present paper is the first one which has results on the
special properties of the bipartite TSP polytope. Our main result is that
a large part of the comb inequalities are not facet defining, although they
are in the case of the complete graph.
%

This paper is organized as follows. The next section contains some
basic notations and notions, also including the description of
the comb inequalities. The main result is discussed
in Section 3. The paper is finished with a short conclusion.

\section{Basic Notations and Notions}
TSP asks for a shortest Hamiltonian tour in an
edge-weighted graph. In its original
form the graph is the complete graph $K_n$. Many valid inequalities of
the TSP polytope are generated by analyzing
how many edges a Hamiltonian tour can have in a special graph
structure. As a matter of fact, even the subtour elimination
constraints can be discussed in that way.
Indeed, for every integer $m$ with $2\leq m<n$, the subtour elimination constraint
 gives the upper bound $m-1$ for the number of edges a tour can have in
 a $K_m$ subgraph of $K_n$. (The formal definition will be given below.)

The basic model of \cite{dfj} uses the binary variables $x_{ij}$
where
\begin{eqnarray*}
x_{ij}=\left\{\begin{array}{ll}1 & \hbox{if the tour goes directly
between cities }i\hbox{ and }j\\ 0 & \hbox{otherwise.}
\end{array}\right.
\end{eqnarray*}
The number of variables $x_{ij}$ is $\left(\begin{array}{c}n\\ 2
\end{array}\right)$
in the undirected case. Therefore the variables $x_{ij}$
and $x_{ji}$ are considered identical for notational convenience.
Each city has two connections. Thus, the variables must satisfy the
constraints
\begin{eqnarray}\label{e07}
\sum_{j=1;\,j\neq i}^n x_{ij}=2.
\end{eqnarray}
In what follows, it is never used that the degrees of vertices must be
exactly~2; we only need that the degrees are at most 2, {\em i.e.}\
\begin{eqnarray}\label{e01}
\sum_{j=1;\,j\neq i}^n x_{ij}\leq 2.
\end{eqnarray}
For all $1\le i\le n$,
inequalities (\ref{e01}) are called {\em degree constraints}. They
still allow the formation of small subtours, which are excluded by the
so-called {\em subtour elimination constraints}:
\begin{eqnarray}\label{e02}
\forall{\cal S}\subset\{1,2,...,n\}, \ 3\leq|{\cal S}|\leq n-1:\;
\sum_{i,j\in{\cal S}}x_{ij}\leq|{\cal S}|-1.
\end{eqnarray}
Only the characteristic vectors of the tours are satisfying both
(\ref{e01}) and (\ref{e02}) among the binary vectors. Continuous
relaxation is obtained if the variables are allowed to take any value
between 0 and 1, {\em i.e.}\ instead of $x_{ij}\in\{0,1\}$
 it is required only that
\begin{eqnarray}\label{e03}
\forall\,1\leq i, j\leq n, \ i\neq j:\;0\leq x_{ij}\leq 1.
\end{eqnarray}
The relaxation is also called Linear Programming (LP) relaxation, especially
if there is a linear objective function, because then the problem becomes an
LP problem. In general, a linear programming problem is an LP relaxation of TSP
if any set of valid inequalities (cuts) of the TSP polytope is a constraint
set of the LP problem. However, the constraints (\ref{e01})--(\ref{e03})
are easy ones, thus it makes no sense to omit any of them.

Let ${\bf x}$ be a feasible solution of the relaxed problem. Then the sum
in (\ref{e02}) is denoted by ${\bf x}({\cal S})$. Thus an
equivalent form of (\ref{e02}) is
$$\forall{\cal S}\subset\{1,2,...,n\}, \ 3\leq|{\cal S}|\leq n-1:\;
{\bf x}({\cal S})\leq|{\cal S}|-1.\eqno(\ref{e02}')$$

The comb inequalities were discovered by Chv\'{a}tal \cite{chvatal2}.
They are the generalizations of Edmonds' 2-matching inequalities. The
{\em comb} is a special graph structure consisting of
two types of sets of vertices: a {\em hand}
and an odd number of {\em teeth}. The intersection of the hand with any
tooth must be non-empty. The graph contains only edges such
that both vertices of the edge are in the same set, {\em i.e.}\ either
in the hand or in one of the teeth. Let $t$ be the number of teeth.
It must be at least 3. In a more formalized way the comb must satisfy
the following constraints. Let ${\cal H}$, ${\cal T}_i$ $(i=1,...,t)$
be the hand and the teeth, respectively. Then
\[{\cal H}\cap{\cal T}_i\neq\emptyset\;\;\mbox{\rm for}\;\;i=1,...,t,\]
\[\forall\,1\leq i,j\leq t, \ i\neq j:\;{\cal T}_i\cap{\cal T}_j=\emptyset,\]
\[t\geq 3 \hbox{  and $t$ is odd.}\]
The comb inequality (see \cite{gutin}, page 59) is
\begin{eqnarray}\label{e04}
{\bf x}({\cal H})+\sum_{i=1}^t {\bf x}({\cal T}_i)\leq |{\cal H}|+
\sum_{i=1}^t |{\cal T}_i|-\frac{3t+1}{2}.
\end{eqnarray}

Notice that each comb has an underlying hypergraph where the hyperedges are the
hand and the teeth.

\section{The Main Results}

Four statements are proven in this section. Lemma 1 is the basic statement. Its
proof shows the technique of all proofs. This lemma is generalized into two
directions in Lemmas 2 and 3. However, Lemmas 2 and 3 are unrelated. A common
generalization of them is Theorem 1 which is the main result of this paper.

All statements are about Hamiltonian circuits in a bipartite graph. It is obvious
that a bipartite graph can have a Hamiltonian circuit only if the sizes of the two
parts are equal. Thus, it may be assumed that the underlying graph is $K_{n,n}$.
However, the assumption that all edges of $K_{n,n}$ are present, is never used in the
proofs.

Let ${\bf x}$ be a feasible solution of the continuous relaxation. It may
have fractional components, of course. Its component $x_{ij}$ belonging to
edge $\{i,j\}$ is called the {\em weight} of the edge. The degree of a vertex
$i$ is the sum of the weights of the edges incident with it. In the bipartite
structure, each vertex belongs to exactly one of two classes. The classes are
denoted by ${\cal C\ell}_1$ and ${\cal C\ell}_2$. Let
${\cal H}^k={\cal H}\cap{\cal C\ell}_k$ for
$k=1,2$.

\subsection{Single Intersections}

\begin{lem}\label{t01} Let $({\cal H},{\cal T}_i\,(i=1,...,t))$ be a comb.
Assume that each vertex of the hand has a separate tooth, {\em i.e.}\
\begin{eqnarray}\label{e08}
|{\cal H}\cap{\cal T}_i|=1,\;\;i=1,...,t
\end{eqnarray}
and
\begin{eqnarray}\label{e09}
{\cal H}\setminus\bigcup_{i=1}^n{\cal T}_i=\emptyset.
\end{eqnarray}
Assume further that a vector ${\bf \bar{x}}$ satisfies the constraints
(\ref{e01})--(\ref{e03}).
Then ${\bf \bar{x}}$ satisfies the comb inequality
(\ref{e04})
in the bipartite case of TSP.

\end{lem}
\proof
Let $p$ and $q$ denote the number of vertices of the hand in classes
${\cal H}^1$ and ${\cal H}^2$, respectively. Thus, $t=p+q$. Let $t=2m+1$
where $m$ is a positive integer.
Without loss of generality it can be assumed that $p\leq m<m+1\leq q$.
We partition the edges of the comb into the following types of sets:
\begin{itemize}
\item Type 1: all edges incident with a vertex $a_i$ of set ${\cal H}^1$,
\item Type 2: all edges of a tooth ${\cal T}_i$ such that they are not incident with
the vertex $a_i={\cal T}_i\cap{\cal H}$ if $a_i\in{\cal H}^1$,
\item Type 3: all edges of a tooth ${\cal T}_i$ with
 ${\cal T}_i\cap{\cal H}=b_i\in{\cal H}^2$.

\end{itemize}
All edges of the comb are elements of exactly one of the sets.
The sets of the three types must satisfy the inequalities as follows.
Type 1 contains all the edges incident with vertex $a_i$. Thus, by the
degree constraints, their
total weight cannot exceed 2:
\begin{eqnarray}\label{e10}
\sum_{u\in{\cal T}_i\setminus\{a_i\}}\bar{x}_{a_i u}+\sum_{u\in{\cal H}^2}
\bar{x}_{a_i u}\leq 2
\end{eqnarray}
where $a_i\in{\cal T}_i$. Type 2 contains all the edges of a subset.
Hence it must satisfy the relevant subtour elimination constraint:
\begin{eqnarray}\label{e11}
\sum_{u,v\in{\cal T}_i\setminus\{a_i\}}\bar{x}_{uv}={\bf \bar{x}}({\cal T}_i
\setminus\{a_i\})\leq \, \mid{\cal T}_i\mid - \, 2.
\end{eqnarray}
Similarly, Type 3 must satisfy a subtour elimination constraint again:
\begin{eqnarray}\label{e12}
\sum_{u,v\in{\cal T}_i}\bar{x}_{uv}={\bf \bar{x}}({\cal T}_i)\leq \,
\mid{\cal T}_i\mid - \, 1
\end{eqnarray}
where ${\cal T}_i$ is the tooth of a vertex in ${\cal H}^2$. The numbers of
inequalities of types (\ref{e10}) and (\ref{e11}) are both $p$. The
number of inequalities of type (\ref{e12}) is $q$. Hence, if all of these
inequalities are summed up then the result is as follows:
\begin{eqnarray}\label{e13}
{\bf \bar{x}}({\cal H})+\sum_{i=1}^t {\bf \bar{x}}({\cal T}_i)
\leq \sum_{i=1}^t |{\cal T}_i| - q\nonumber\\
=|{\cal H}|+\sum_{i=1}^t |{\cal T}_i| - p - 2q\\
\leq |{\cal H}|+\sum_{i=1}^t
|{\cal T}_i|-
\frac{3t+1}{2}.\nonumber
\end{eqnarray}
\qed

According to the definition of the comb, the hand may contain vertices without
tooth. Lemma \ref{t01} can be extended to this case as well.

\begin{lem}\label{t03} Let $({\cal H},{\cal T}_i\,(i=1,...,t))$ be a comb.
Assume that no tooth contains two or more vertices of the hand, {\em i.e.}
\begin{eqnarray}\label{e14}
|{\cal H}\cap{\cal T}_i|= 1,\;\;i=1,...,t.
\end{eqnarray}
Assume further that a vector ${\bf \bar{x}}$ satisfies the constraints
(\ref{e01})--(\ref{e03}).
Then ${\bf \bar{x}}$ satisfies the comb inequality
in the bipartite case of TSP.
\end{lem}
\proof
Let $p$ and $q$ be the number of the elements of ${\cal H}^1$ and ${\cal H}^2$, respectively.
This time it is not assumed that $p<q$.
Assume that ${\cal G}^1$ and ${\cal G}^2$ are the subsets of ${\cal H}^1$ and ${\cal H}^2$ such that they have tooth.
Let the numbers of the elements of the sets ${\cal G}^1$ and ${\cal G}^2$
be $p'$ and $q'$, respectively; then $t=p'+q'$.

A method similar to the proof of the previous lemma is applied.
In what follows, we denote a vertex of ${\cal G}^1\cap{\cal T}_i$ by $a_i$.
All edges of the comb structure belong to one of the
following four disjoint sets:
\begin{itemize}
\item Type 1: all edges incident with a vertex $a_i$ of set ${\cal G}^1$,
\item Type 2: all edges incident with a vertex $b_i\in{\cal H}^1\setminus{\cal G}^1$,
\item Type 3: all edges of a tooth ${\cal T}_i$ such that they are not incident with
the vertex $a_i={\cal T}_i\cap{\cal H}$ if $a_i\in{\cal H}^1$,
\item Type 4: all edges of a tooth ${\cal T}_i$ with
 ${\cal T}_i\cap{\cal H}=c_i\in{\cal H}^2$.
\end{itemize}

The same
inequality as (\ref{e10}) can be claimed for the elements of ${\cal G}^1$ as follows:
\begin{eqnarray}\label{e15}
\sum_{u\in{\cal T}_i\setminus\{a_i\}}\bar{x}_{a_i u}+\sum_{u\in{\cal H}^2}
\bar{x}_{a_i u}\leq 2.
\end{eqnarray}
The number of these inequalities is $p'$. The degree constraints of
the elements of ${\cal H}^1\setminus {\cal G}^1$ are simpler:
\begin{eqnarray}\label{e16}
\sum_{u\in{\cal H}^2} \bar{x}_{b_i u}\leq 2.
\end{eqnarray}
The number of these inequalities is $p-p'$. The subtour elimination constraints
of the teeth belonging to ${\cal G}^1$ vertices are exactly
the same as the inequalities (\ref{e11}):
\begin{eqnarray}\label{e17}
\sum_{u,v\in{\cal T}_i\setminus\{a_i\}}\bar{x}_{uv}={\bf \bar{x}}({\cal T}_i
\setminus\{a_i\})\leq \, \mid{\cal T}_i\mid - \, 2.
\end{eqnarray}
Their number is again $p'$. Finally, the subtour elimination constraints of the teeth belonging to ${\cal G}^2$ vertices are as follows:
\begin{eqnarray}\label{e18}
\sum_{u,v\in{\cal T}_i}\bar{x}_{uv}={\bf \bar{x}}({\cal T}_i)\leq \,
\mid{\cal T}_i\mid - \, 1.
\end{eqnarray}
Their number is $q'$. If the inequalities (\ref{e15})-(\ref{e18}) are summed up then the inequality
\begin{eqnarray}\label{e19}
{\bf \bar{x}}({\cal H})+\sum_{i=1}^t {\bf \bar{x}}({\cal T}_i)
\leq \sum_{i=1}^t |{\cal T}_i| - q'+2(p-p')
\end{eqnarray}
is obtained. Notice that $|{\cal H}|=p+q$. Thus the right-hand side of
(\ref{e19}) can be reformulated as follows:
\begin{eqnarray*}
\sum_{i=1}^t |{\cal T}_i| - q'+2(p-p')=|{\cal H}|+\sum_{i=1}^t |{\cal T}_i| - q'-2p'+p-q.
\end{eqnarray*}
If
\begin{eqnarray}\label{e20}
\frac{3t+1}{2}\leq q'+2p'-p+q
\end{eqnarray}
then the statement holds. Assume that the opposite of
(\ref{e20}) is true, which is equivalent to
\begin{eqnarray}\label{e21}
3p'+3q'+1>2q'+4p'-2p+2q.
\end{eqnarray}
The procedure is repeated such that the roles of ${\cal H}^1$ and
${\cal H}^2$ are interchanged. Then either
\begin{eqnarray*}
\frac{3t+1}{2}\leq p'+2q'-q+p
\end{eqnarray*}
and the statement is true, or
\begin{eqnarray}\label{e22}
3p'+3q'+1> 2p'+4q'-2q+2p.
\end{eqnarray}
By summing up the inequalities (\ref{e21}) and (\ref{e22}) the relation \begin{eqnarray}\label{e23}
6p'+6q'+2 > 6p'+6q'.
\end{eqnarray}
is obtained. Now (\ref{e23}) shows that the difference of the left-hand side and the right-hand side
in both (\ref{e21}) and (\ref{e22}) is 1. However it is not possible
as there are even numbers on both sides of both inequalities.
\qed

\subsection{Multiple Intersections}

Notice that in the previous results,
 one of the classes has vertices in less than the half of the
 intersections of the teeth and the hand, as each intersection consists
of a single vertex.

In the next step it is allowed that the intersection of the hand and
a tooth have more than one element. However, (\ref{e09}) is still required, {\em i.e.}
each element of the hand is in a tooth as well. Further on, vertices of one of
the classes occur in less than the half of the intersections of teeth and hand.

\begin{lem}\label{t02}
Let $({\cal H},{\cal T}_i\,(i=1,...,t))$ be a comb.
Let $p$ and $q$ be two positive integers such that $t=p+q$ with $p<q$.
Assume that ${\cal H}^1\cap{\cal T}_i\neq\emptyset$ for $i=1,...,p$
and  ${\cal H}^1\cap{\cal T}_i=\emptyset$ for $i=p+1,...,t$.
Assume further that Constraint (\ref{e09}) holds.
If a vector ${\bf \bar{x}}$ satisfies the constraints
(\ref{e01})--(\ref{e03})
then ${\bf \bar{x}}$ satisfies the comb inequality
in the bipartite case of TSP.
\end{lem}
\proof
We introduce the following notation for the intersections
 $\mid {\cal H}^1\cap{\cal T}_i\mid$ and $\mid {\cal H}^2\cap{\cal T}_i\mid$,
 which will be used also in the more general result later:
\begin{itemize}
 \item $\mid {\cal H}^1\cap{\cal T}_i\mid=1+s_i\;(i=1,...,p)$,
 \item $\mid {\cal H}^2\cap{\cal T}_i\mid=r_i\;(i=1,...,p)$,
 \item $\mid {\cal H}^2\cap{\cal T}_i\mid=1+r_i\;(i=p+1,...,t)$,
 \item $\mid {\cal H}^1\cap{\cal T}_i\mid=0\;(i=p+1,...,t)$.
\end{itemize}
In accordance with the conditions of the assertion, it
 is assumed here that $s_i$ and $r_i$ are nonnegative integers;
 the last row is recalled for the sake of completeness.

The method of Lemma \ref{t01} is applied again. However, there is one significant
difference. In the case of Lemma \ref{t01} each intersection of the hand and a tooth
consists of a single vertex, thus, every edge on the left-hand side of (\ref{e04})
appears only once. There are edges in ${\cal T}_i\cap{\cal H}$ if $1\leq i\leq p$
and $r_i>0$. These edges occur on the left-hand side of (\ref{e04}) twice, once in
${\bf x}({\cal H})$ and once in ${\bf x}({\cal T}_i)$. Thus, they must be covered
twice. Now we define the types of the edges of the comb as follows:
\begin{itemize}
\item Type A: edges in ${\cal T}_i\setminus{\cal H}$, where $1\leq i\leq p$,
\item Type B: edges with one vertex in ${\cal T}_i\setminus{\cal H}$ and
with another vertex in ${\cal T}_i\cap{\cal H}^1$, where $1\leq i\leq p$,
\item Type C: edges in ${\cal T}_i\cap{\cal H}$,
where $1\leq i\leq p$,
\item Type D: edges with one vertex in ${\cal T}_i\setminus{\cal H}$ and
with another vertex in ${\cal T}_i\cap{\cal H}^2$, where $1\leq i\leq p$,
\item Type E: edges in the hand connecting two different teeth,
\item Type F: edges with one vertex in ${\cal T}_i\setminus{\cal H}$ and
with another vertex in ${\cal T}_i\cap{\cal H}^2$, where $p+1\leq i\leq t$,
\item Type G: edges in ${\cal T}_i\setminus{\cal H}$, where $p+1\leq i\leq t$.
\end{itemize}
Type C edges must be covered twice.
The types of the subsets covering the edges of the comb are as follows:
\begin{itemize}
\item Type 1: The edges incident with one vertex in
 ${\cal T}_i\cap{\cal H}^1$, where $1\leq i\leq p$. These
 sets of edges cover edge Types B, C, and E.
\item Type 2: The edges incident with one vertex in
${\cal T}_i\cap{\cal H}^2$, where $1\leq i\leq p$. These
sets of edges cover edge Types C, and D
 (and possibly may contain some edges of Type E,
 which we will disregard in the computation without loss of generality).
\item Type 3: The edges of ${\cal T}_i\setminus{\cal H}$, where $1\leq i\leq p$.
These sets cover Type A edges.
\item Type 4: The edges of ${\cal T}_i$, where $p+1\leq i \leq t$.
These sets cover Types F, and G edges.
\end{itemize}

Let $a$ be a Class 1 vertex in the intersection of ${\cal H}$
and ${\cal T}_i$. The form of the inequalities of Type 1 set is
based on the degree constraint. It is
\begin{eqnarray}\label{e24}
\sum_{u\in{\cal T}_i} x_{au}+\sum_{u\in{\cal H}^2\setminus{\cal T}_i}x_{au}\leq 2.
\end{eqnarray}
The number of these constraints is $p+\sum_{i=1}^p s_i$. There are
$\sum_{i=1}^p r_i$ inequalities of Type 2 and their form is
\begin{eqnarray}\label{e25}
\sum_{u\in{\cal T}_i}
x_{au}+\sum_{u\in{\cal H}^1\setminus{\cal T}_i}x_{au}\leq 2,
\end{eqnarray}
where $a\in{\cal T}_i\cap{\cal H}^2$ with $1\leq i\leq p$.
The number of Type 3 constraints is $p$. The related
inequalities are based on the subtour elimination constraints as follows:
\begin{eqnarray}\label{e26}
\sum_{u,v\in{\cal T}_i\setminus{\cal H}}
x_{uv}\leq \mid{\cal T}_i\mid-s_i-r_i-2.
\end{eqnarray}
here and also in 2nd line of text
 below, and also in next display \\
The right-hand side is obtained from the fact that number of elements of
${\cal T}_i\setminus{\cal H}$ is $\mid{\cal T}_i\mid-s_i-r_i-1$ if
$1\leq i\leq p$.
Finally, Type 4 constraints are not changed. Their number is $q$. The related
inequalities are
\begin{eqnarray}\label{e27}
\sum_{u,v\in{\cal T}_i}
x_{uv}\leq \mid{\cal T}_i\mid-1,
\end{eqnarray}
where $p+1\leq i\leq t$. By summing up all these inequalities, a new constraint is
obtained as follows:
\begin{eqnarray*}\label{e28}
{\bf x}({\cal H})+\sum_{i=1}^t {\bf x}({\cal T}_i)
\;\;\;\\ \leq 2p+2\sum_{i=1}^{p}s_i
+2\sum_{i=1}^{p}r_i+\sum_{i=1}^t |{\cal T}_i|-\sum_{i=1}^{p}s_i
-\sum_{i=1}^{p}r_i -2p-q\;\;\;\\ =\sum_{i=1}^t |{\cal T}_i|+\sum_{i=1}^{p}s_i
+\sum_{i=1}^{p}r_i -q\;\;\;\\ =\sum_{i=1}^t |{\cal T}_i|+\sum_{i=1}^{p}s_i
+\sum_{i=1}^{p}r_i -q+\mid{\cal H}\mid-p-\sum_{i=1}^{p}s_i
-\sum_{i=1}^{p}r_i -q-\sum_{i=p+1}^t r_i\;\;\;\\
=\mid{\cal H}\mid + \sum_{i=1}^t |{\cal T}_i|-p-2q-\sum_{i=p+1}^t r_i\;\;\;\\
\leq\mid{\cal H}\mid + \sum_{i=1}^t |{\cal T}_i|-p-2q\;\;\;\\
\leq\mid{\cal H}\mid + \sum_{i=1}^t |{\cal T}_i|-\frac{3t+1}{2}.
\end{eqnarray*}
Thus, the comb inequality holds.
\qed

The following example shows that if vertices from both classes occur in more than the half
of the teeth, then the statement is not true in general.
The comb is described in Table 1.

\begin{center}
\begin{tabular}{|c|c|c|}
  \hline
   & ${\cal C\ell}_1$ & ${\cal C\ell}_2$ \\ \hline
  ${\cal H}$ & $a,b$ & $c,d$ \\
  ${\cal T}_1$ & $a$ & $e$ \\
  ${\cal T}_2$ & $b,g$ & $c,f$ \\
  ${\cal T}_3$ & $h$ & $d$ \\
  \hline
\end{tabular}
\centerline{{\bf Table 1}}
\end{center}

The weights of the edges, {\em i.e.}\ the values of the variables, are
as follows: $1=x_{ae}=x_{bc}=x_{fg}=x_{dh}$, $0.5=x_{ac}=x_{ad}=x_{bd}=x_{bf}=x_{cg}$.
Hence, ${\bf x}({\cal H})=2.5$,  ${\bf x}({\cal T}_1)=1$, ${\bf x}({\cal T}_2)=3$,
${\bf x}({\cal T}_3)=1$. Thus, the value of the left-hand side of the comb inequality is
7.5. The value of the right-hand side is $4+(2+4+2)-5=7$, where $5=\frac{3t+1}{2}$.

In the next step the hand may contain vertices without tooth.
We use the notation introduced above for the pattern of intersections
of the teeth with the hand. We recall it for easier reference in the next table:

\begin{center}
\begin{tabular}{|c||c|c|}
  \hline
   & ${\cal H}^1$ & ${\cal H}^2$ \\ \hline\hline
  $1\le i\le p$ & $1+s_i$ & $r_i$ \\ \hline
  $p+1\le i\le t$ & $0$ & $1+r_i$ \\
  \hline
\end{tabular}
\end{center}

\begin{thm}\label{t04}
Let $({\cal H},{\cal T}_i\,(i=1,...,t))$ be a comb with intersection
pattern defined above.
Let $p$ and $q$ be two positive integers such that $t=p+q$.
Assume that $w$ and $y$ are the number of the elements of ${\cal H}^1$
 and ${\cal H}^2$ such that they have no tooth. If
\begin{eqnarray}\label{e29}
w\leq
 y + \frac{q-(p+1)}{2}
+\sum_{i=p+1}^t r_i
\end{eqnarray}
then the comb inequality is satisfied for all
vectors ${\bf \bar{x}}$ satisfying the constraints
(\ref{e01})--(\ref{e03}).

\end{thm}
\proof
We use essentially the same edge types and subset types as in the proof of Lemma \ref{t02},
with two exceptions. In the current case, we redefine
\begin{itemize}
\item Type E: edges in the hand connecting two different teeth, and also those edges
 inside ${\cal H}$ which are incident with at least one vertex without tooth.
\end{itemize}
Similarly, the edge set of Type 1 should contain the additional edges of Type~E. Thus

\begin{itemize}
\item Type 1: The edges incident with one vertex in
 ${\cal T}_i\cap{\cal H}^1$, where $1\leq i\leq p$, and the edges incident with the $w$ elements of ${\cal H}^1$ without tooth.
\end{itemize}

The form of the inequalities of the Type 1 set is based on the degree constraint, and
compared to the previous proof it is complemented with those for vertices in
 ${\cal H}^1\setminus(\bigcup{\cal T}_i)$.
 The number of the constraints in the form of (\ref{e24}) is $p+\sum_{i=1}^p s_i$.
  Moreover, if $a$ is an element of ${\cal H}^1$ without tooth, then the form of the inequalities of Type 1 is
\begin{eqnarray}\label{e30}
\sum_{u\in{\cal H}^2}x_{au}\leq 2.
\end{eqnarray}
The number of these constraints is $w$.

The inequalities of all other types are the same as in the case of Lemma~\ref{t02}.

By summing up all these inequalities, a new constraint is
obtained as follows:
\begin{eqnarray*}\label{e31}
{\bf x}({\cal H})+\sum_{i=1}^t {\bf x}({\cal T}_i)\;\;\;\\
\leq 2p+2\sum_{i=1}^{p}s_i+2\sum_{i=1}^{p}r_i+2w+\sum_{i=1}^t |{\cal T}_i|-\sum_{i=1}^{p}s_i-\sum_{i=1}^{p}r_i -2p-q\;\;\;\\
=\sum_{i=1}^t |{\cal T}_i|+\sum_{i=1}^{p}s_i+\sum_{i=1}^{p}r_i -q+2w\;\;\;\\
=\sum_{i=1}^t |{\cal T}_i|+\sum_{i=1}^{p}s_i+\sum_{i=1}^{p}r_i -q+2w\;\;\;\\
+\mid{\cal H}\mid-p-\sum_{i=1}^{p}s_i-\sum_{i=1}^{p}r_i -q-\sum_{i=p+1}^t r_i-w-y\;\;\;\\
=\mid{\cal H}\mid + \sum_{i=1}^t |{\cal T}_i|-p-2q+w-y-\sum_{i=p+1}^t r_i.\;\;\;
\end{eqnarray*}
Thus, the comb inequality holds whenever
\begin{eqnarray}\label{e32}
-p-2q+w-y-\sum_{i=p+1}^t r_i\leq-\frac{3t+1}{2},
\end{eqnarray}
and it means that
\begin{eqnarray}\label{e33}
w\leq\frac{q}{2}-\frac{p}{2}-\frac{1}{2}+y+\sum_{i=p+1}^t r_i
\end{eqnarray}
which is equivalent to (\ref{e29}).\qed


The next example shows that if the inequality (\ref{e33}) does not hold,
 then the statement can be violated.
The comb is described in Table 2.

\begin{center}
\begin{tabular}{|c|c|c|}
  \hline
   & ${\cal C\ell}_1$ & ${\cal C\ell}_2$ \\ \hline
  ${\cal H}$ & $a,b$ & $c,d,e$ \\
  ${\cal T}_1$ & $a,f$ & $c,i$ \\
  ${\cal T}_2$ & $g$ & $d$ \\
  ${\cal T}_3$ & $h$ & $e$ \\
  \hline
\end{tabular}
\centerline{{\bf Table 2}}
\end{center}

Without the vertex $e$, this example is covered in Lemma \ref{t02}, because $p=1$, $q=2$.

In this case $e$ has no tooth, $w=1$ and $y=0$, $p=1$, $q=2$.

The weights of the edges, {\em i.e.} the values of the variables, are
as follows: $1=x_{ac}=x_{fi}=x_{dg}=x_{eh}$, $0.5=x_{ad}=x_{ai}=x_{bc}=x_{bd}=x_{cf}$.
Hence, ${\bf x}({\cal H})=3.5$,  ${\bf x}({\cal T}_1)=3$, ${\bf x}({\cal T}_2)=1$,
${\bf x}({\cal T}_3)=1$. Thus, the value of the left-hand side of the comb inequality is
8.5. The value of the right-hand side is $5+(4+2+2)-5=8$, where $5=\frac{3t+1}{2}$.


\subsection{Symmetric Case and Hierarchy}


%

If the notations of Theorem 1 are applied to Lemma 3 then $w=y=0$.
Thus, the theorem is more general than Lemma 3.
It is even more obvious that Lemma 1 follows from Lemma 2.
On the other hand, Lemmas 2, and 3 are unrelated.

In this short subsection we consider the situation where each tooth
intersects the hand in only one vertex class. We prove that
the side condition (\ref{e28}) can be omitted in this case.
This will immediately imply that the theorem is a generalization
of Lemma 2 as well. Hence, Theorem 1 is the most general statement
proven in this paper.

Compared to Theorem 1, here we consider the following simplified
intersection pattern:

\begin{center}
\begin{tabular}{|c||c|c|}
  \hline
   & ${\cal H}^1$ & ${\cal H}^2$ \\ \hline\hline
  $1\le i\le p$ & $1+s_i$ & $0$ \\ \hline
  $p+1\le i\le t$ & $0$ & $1+r_i$ \\
  \hline
\end{tabular}
\end{center}

\begin{thm}\label{t05}
Let $({\cal H},{\cal T}_i\,(i=1,...,t))$ be a comb with intersection
pattern defined above, \emph{i.e.}\ each tooth
intersects the hand in only one vertex class.
Let $p$ and $q$ be two positive integers such that $t=p+q$.
Then the comb inequality holds for all
vectors ${\bf \bar{x}}$ satisfying the constraints
(\ref{e01})--(\ref{e03}).

\end{thm}
\proof
Assume that $w$ and $y$ are the number of the elements of ${\cal H}^1$
and ${\cal H}^2$ such that they have no tooth. The proof of Theorem 1
implies the current assertion unless (\ref{e32}) is violated. Since
$t$ is odd, and now we have $\sum_{i=1}^p r_i=0$,
violation of (\ref{e32}) is equivalent to
\begin{eqnarray}\label{e34}
\frac{3t+1}{2} \geq p+2q-w+y+1.
\end{eqnarray}
Due to the symmetry of the intersection pattern, we may also switch the role
of the two vertex classes, and analogously derive
\begin{eqnarray}\label{e35}
\frac{3t+1}{2} \geq 2p+q+w-y+1.
\end{eqnarray}
But then the sum of (\ref{e34}) and (\ref{e35}) yields
\begin{eqnarray}
3t+1 \geq 3p+3q+2 = 3t+2,
\end{eqnarray}
a contradiction that completes the proof of the theorem.
\qed


\section{Conclusions}
The main result of this study is that many comb inequalities are not facet defining
for the polytope of the symmetric TSP in the case of bipartite graphs. It is not
a contradiction to the results of \cite{clique1} and \cite{n02} as the
original theorem on the facet defining property concerns the TSP
polytope of the complete graph $K_n$. Thus the structure of the TSP
polytope depends on the underlying graph structure.

\end{sloppypar}
\end{document}